  \documentclass[5p]{elsarticle}

\usepackage{amssymb,amsmath}

\def\LM#1#2{\left|\begin{array}{l}{#1}\\[1ex]{#2}\end{array}\right.}

\begin{document}

\begin{frontmatter}

\title{Fast, Accurate and Robust Adaptive Finite Difference Methods for Fractional Diffusion Equations: The Size of the Timesteps does Matter}

\author{ Santos B. Yuste  and         Joaqu\'{\i}n Quintana-Murillo }
\address{ Departamento de F\'{\i}sica \\
        Universidad de Extremadura, E06071 Badajoz, Spain}

\begin{abstract}
 The computation time required by standard finite difference methods  with fixed timesteps  for solving  fractional diffusion equations is usually very large because the number of operations required to find the solution scales as the square of the number of timesteps.  Besides, the solutions of these problems usually involve markedly different time scales, which leads to quite inhomogeneous numerical errors.  A natural way to address these difficulties is by resorting to adaptive numerical methods where the size of the timesteps is chosen according to the behaviour of the solution.  A key feature of these methods is then the efficiency of   the adaptive algorithm employed to \emph{dynamically} set the size of every timestep. Here we discuss two adaptive methods based on the step-doubling technique. These  methods are, in many cases, immensely faster  than the corresponding standard method with fixed timesteps and they allow  a tolerance level to be set for the numerical errors that turns out to be a good indicator of the actual errors.
\end{abstract}


\end{frontmatter}


\section{Introduction}

Interest in fractional calculus is flourishing,  to a large extent, due to its usefulness  as a mathematical tool for tackling an increasing variety of scientific problems usually associated with complex systems that show  some kind of long-term memory. In  control engineering, fractional calculus have been successfully employed for many years.  In statistical physics, fractional calculus is  extremely useful in the study of some systems in which the diffusive processes are anomalous.  In particular, it can be rigorously proved that  fractional diffusion equations are the appropriate equations to describe the spread of some classes of continuous time random walkers  in the same way that  normal  diffusion equations  describe the diffusion of  Brownian  walkers (see, for example, the review chapter of Metzler and Jeon in Ref.~\cite{Klafter2011a}).

Of course, the utility of  fractional calculus is linked to our ability to extract useful information about the systems concerned from this formalism. Fortunately, many of the long-time well-established analytical methods employed to study normal diffusion equations can be adapted to fractional diffusion equations. For example,  in many cases fractional and normal diffusion equations can be solved similarly in Fourier-Laplace space.  However, as is the case too with  non-fractional problems, it is also very convenient (or even indispensable)  to have at our disposal  numerical procedures for studying these equations,  and thus get information about the systems they describe.   The study of numerical methods for solving fractional equations has been an area  of quite active research in recent years (see   relatively recent accounts of the literature on this issue in Refs. \cite{YangThesis2010} and \cite{Li2012}). As also for non-fractional differential equations, finite difference methods are one of the most important classes of  numerical methods for solving fractional partial differential equations.

Usually, finite difference methods developed for fractional diffusion equations employ   uniform time discretization, i.e., fixed timesteps \cite{Li2012}.  But methods of this kind have two main drawbacks:  they are slow and their accuracy is inconsistent. In fact,  they become \emph{increasingly} slower as time goes by: the CPU time required to get the solution at time $t$  grows as the \emph{square} of $t$ (i.e., the arithmetic complexity of these algorithms is of the order $\Delta^{-2}$, with $\Delta$ being the size of the timestep).
This difficulty has been recognized for a long time, and some procedures have been proposed to alleviate it.  The most obvious is to increase the order of accuracy of the numerical method so that larger timesteps can be used without losing the accuracy of the solutions \cite{Li2012}.  Another approach is based on the so-called ``short memory principle'' \cite{Deng2007a}  that, in summary, either tries to cap the number of required operations per step assuming that the influence of the previous values of the solution for times far from the present time can be neglected \cite{Podlubny1999}, or takes advantage of the way in which the kernel of the fractional derivative decays to get arithmetic complexity of order $\Delta^{-1} \log \Delta^{-1}$ \cite{Ford2001}.
Another problem of standard methods with fixed timesteps, one that is rarely noticed, is that the accuracy of their numerical results changes  strongly (even by orders of magnitude, see the figures in Section \ref{AdapMet}) over the time interval of integration.  The reason for this behaviour can be traced back to the typical behaviour of the solutions of many fractional diffusion equations. In many cases, these solutions can be written as a superposition of generalized Fourier modes that decay as Mittag-Leffler functions. But, as is well known \cite{Mainardi2000}, these functions decay very fast for short times and very slowly for longer times. This behaviour is, in many cases, inherited by the full solution. Indeed, in these cases, to use fixed timesteps to deal with such different time regimes seems a poor choice.

These problems regarding the speed and accuracy of finite difference methods can be alleviated by using methods with adaptive timesteps. This kind of method has the great advantage that the size of the timesteps can be chosen  according to the behaviour of the solution. Ideally, a good adaptive method, as Press et al. say for ODE integrators,  ``should exert some adaptive control over its own progress, making frequent changes in its stepsize \ldots Many small steps should tiptoe through treacherous terrain, while a few great strides should speed through smooth uninteresting countryside \ldots[so as to] achieve some predetermined accuracy in the solution with minimum computational effort'' \cite{Press:2007:NRE:1403886}.
In order to construct this kind of method for fractional diffusion equations, two key ingredients are required:  first,   a finite difference method that can work with variable timesteps, and second,   a procedure for choosing the size of the timesteps.
Finite difference methods that can work with variable timesteps are scarce. Some examples are the matrix approach on non-equidistant grids by Podlubny et al. \cite{Podlubny2013}, a  generalized Crank–-Nicolson method   by Mustapha et al. \cite{Mustapha2010,Mustapha2012}, and a non-uniform L1 time discretization  \cite{Yuste2012,Quintana-Murillo2013,Zhang2014}.
The finite difference method we employ in this paper is an unconditionally stable implicit method discussed in Ref.  \cite{Yuste2012}; the adaptive control procedure is based in the so-called  step-doubling technique \cite{Press:2007:NRE:1403886}.

The paper is organized as follows.  In Section \ref{AlgoVar} we present an unconditionally stable finite difference scheme that is able to solve  fractional diffusion equations by employing  variable timesteps. In Section \ref{AdapMet}  we present two adaptive algorithms for choosing the size of the timesteps and we analyze their speed and accuracy. In Section \ref{twoExamples}  we provide two examples that show some relevant features and strengths of our adaptive methods.  We end with some remarks and conclusions.

\section{Algorithm with non-uniform timesteps}
\label{AlgoVar}

The equation we consider is a one-dimensional fractional diffusion equation in the Caputo form
\begin{equation}
 \frac{\partial^\gamma u}{\partial t^\gamma}   =
  K \frac{\partial^2 u}{\partial x^2}   +f(x,t)
 \label{ecu1}
\end{equation}
where $f(x,t)$ is a source term  and
\begin{equation}
\label{CaputoRL2b}
  \frac{\partial^\gamma}{\partial t^\gamma} y(t)\equiv
\frac{1}{\Gamma(1-\gamma)} \int_0^t d\tau  \frac{1}{(t-\tau)^{\gamma}}
\frac{dy(\tau)}{d\tau},\quad 0<\gamma<1,
\end{equation}
is the Caputo fractional derivative \cite{Podlubny1999}. The extension of our procedure to other spatial dimensions and to other  equations with  terms involving standard non-fractional  spatial derivatives (e.g., the fractional Fokker-Planck equation \cite{Metzler2000}) is straightforward.

In this paper,  the Caputo time derivative is discretized by means of a direct  generalization of the well-known fractional L1 formula \cite{OldhamSpanier1974} to the case of non-uniform meshes \cite{Yuste2012}.  This non-uniform time discretization is a key part of our approach. Because our purpose is to study adaptive methods that tackle the difficulties associated with the fractional nature of the time-derivative operator by employing non-uniform timesteps, we limit ourselves to the simple  discretization of the non-fractional part of the equation (the Laplacian operator) by means of the three-point centred formula. For the case of uniform timesteps  the present method becomes the numerical scheme discussed by Liu et al. \cite{Liu2006} and Murio \cite{Murio2008}. It can be proved that the method is unconditionally stable \emph{regardless} of the size of the (non-uniform) timesteps employed \cite{Yuste2012}. A key aspect of the present method is the way in which the fractional derivative is discretized on a non-uniform temporal mesh; the discretization of the \emph{non}-fractional spatial operator can be implemented straightforwardly by means of standard procedures of non-fractional finite difference methods \cite{LeVeque2007}.
For the sake of completeness, we shall give here the main formulas of this finite difference scheme on  non-uniform temporal meshes.

Let $(x_j,t_m)$ be  the coordinates of the $(j,m)$ node of the mesh of the space-time region where one wants to obtain the numerical solution of the fractional equation. We will denote by  $U_j^{(m)}$  the numerical estimate provided by the difference methods of the exact solution $u(x_j,t_m)=u_j^{(m)}$. Next, we replace the continuous operators of the fractional equation by suitably chosen difference operators:
 \begin{align}
 \frac{\partial^\gamma}{\partial t^\gamma} u(x,t_{n})=&\frac{1}{\Gamma(2-\gamma)} \sum_{m=0}^{n-1} T_{m,n}^{(\gamma)} \left[ u(x,t_{m+1})-u(x,t_{m})\right] \nonumber\\
  &+R_{t_n}(x)
\end{align}
where \cite{Yuste2012}
 \begin{align}
T_{m,n}^{(\gamma)}   &=\frac{\int_{t_m}^{t_{m+1}} (t_n-\tau)^{-\gamma}\,d\tau}{t_{m+1}-t_{m}}   \nonumber\\
 &=     \frac{(t_{n}-t_{m})^{1-\gamma} - (t_{n}-t_{m+1})^{1-\gamma}}{t_{m+1}-t_{m}}, \quad m\le n-1.
\label{Tmn}
\end{align}
The truncation error $R_{t_n}(x)$ is bounded by a quantity of order $t_n^{1-\gamma} \Delta_\text{max}$ with   $\displaystyle{\Delta_\text{max}=\max_{0\le m\le n-1}\left(t_{m+1}-t_{m} \right)}$. Recently, Zhang et al.~\cite{Zhang2014} have reported an improved bound that is the sum of one term of order $\Delta_n^{2-\gamma}$ and another of order $\Delta_\text{max}^{2}\Delta_n^{-\gamma}$.
 The Laplacian operator is given by the three-point centred formula:
\begin{equation}
\label{Lapladiscre2}
  \frac{\partial^2}{\partial x^2}  u(x_j,t)=\frac{u(x_{j+1},t)-2u(x_{j},t)+u(x_{j+1},t)}{(\Delta x)^2} + R_{x_j}(t).
\end{equation}
The truncation error  $R_{x_j}(t)$ is of order $(\Delta x)^2 $.  Neglecting the truncation errors and multiplying the equation by $(t_{n} -t_{n-1})^\gamma$ one gets the difference approximation to the continuous equation \cite{Yuste2012}:
   \begin{equation}
   \label{metImp1}
   \sum_{m=0}^{n-1} \tilde T_{m,n}^{(\gamma)} \left[ U^{(m+1)}_{j} -U^{(m)}_{j}\right]
     =S_n  [U^{(n)}_{j+1}- 2 U^{(n)}_{j} +U^{(n)}_{j-1}]+  F^{(n)}_{j}
\end{equation}
where \cite{Yuste2012,Yuste2014Corrigendum}
\begin{align}
S_n&= \Gamma(2-\gamma)  K \frac{ \left(t_{n} -t_{n-1}\right)^\gamma}{(\Delta x)^2},\\
\label{tilTdef}
\tilde T_{m,n}^{(\gamma)}&=  (t_{n}-t_{n-1})^{\gamma}T_{m,n}^{(\gamma)}\, ,\\
F^{(n)}_{j}&=\Gamma(2-\gamma)(t_{n}-t_{n-1})^{\gamma}f(x_n,t_n).
\end{align}
Reordering  (\ref{metImp1}) one gets the following (implicit) finite  difference scheme:
\begin{align}
  &-S_n \,  U^{(n)}_{j+1} +(1+2S_n) U^{(n)}_{j}-S_n \, U^{(n)}_{j-1} \nonumber \\
  &= U^{(n-1)}_{j} -   \sum_{m=0}^{n-2} \tilde T_{m,n}^{(\gamma)} \left[ U^{(m+1)}_{j} -U^{(m)}_{j}\right]
  +F^{(n)}_{j},
  \label{ImpliAlg1}
\end{align}
which can be written in matrix-vector form as
\begin{equation}
A \, U^{(n)}=G\left(U^{(n-1)},U^{(n-2)},\ldots,U^{(0)},F^{(n)},t_n\right),
\end{equation}
where $U^{(m)}$ stands for the vector $\{U_j^{(m)}\}$.  This equation, $A U=G$,  is formally identical to that of the non-fractional differential equation, and its solution $U=A^{-1} G$ can be obtained very efficiently by means of the Thomas algorithm  because $A$  is  a tridiagonal matrix. The key difference with respect to non-fractional algorithms is that, in order to evaluate $G$,  the numerical solution $U^{(m)}$ for \emph{all} the $n$ previous time values has to be employed [see the right-hand side of \eqref{ImpliAlg1}], while for non-fractional equations (i.e., for $\gamma=1$) only the solution at the previous value $U^{(n-1)}$ is required.  This explains why finite difference methods are increasingly slow: the computational cost of going from the solution at time  $t_{n-1}$ to the solution at time $t_n$ grows as $n$, i.e., as  the number of terms of the sum that defines $G$, which implies that the computational cost for going from $t_0$ to $t_n$ grows as $n^2$.

\section{Adaptive methods}
\label{AdapMet}

\begin{figure}[t]
\begin{center}
\includegraphics[width=0.45\textwidth,angle=0]{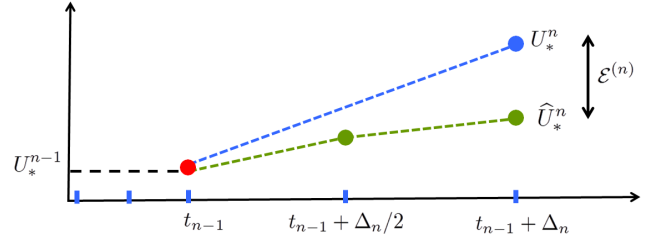}    
\end{center}
\caption{\label{fig:diferenciaU}
Scheme of the step-doubling technique.  The solution at time $t_n$ is obtained by means of (i) a full timestep of size $t_n-t_{n-1}$ and (ii) by means of two steps of size $(t_n-t_{n-1})/2$. The difference $\mathcal{E}^{(n)}$ between both solutions is used as an indicator of the numerical error.
}
\end{figure}

In the previous section we have presented a finite difference method that can work with variable timesteps. This is the first key ingredient of our adaptive method.
The second  ingredient  is a procedure for choosing the size of the timesteps according to the behaviour of the solution.  In this paper, we shall  consider two methods: the trial and error (T\&E) step-doubling algorithm \cite{Quintana-Murillo2013}, and the predictive  step-doubling algorithm. Both algorithms are based on the step doubling technique \cite{Press:2007:NRE:1403886}:  the numerical solution at a given time $t_{n}$ is evaluated twice, first  employing a full step $\Delta_n=t_{n}-t_{n-1}$ and, next, independently, using two half steps of size $\Delta_n/2$;   the difference  $\mathcal{E}^{(n)}$  between the two numerical estimates of the solution,  ${U}_k^{(n)}$ and $\widehat{U}_k^{(n)}$, respectively (see Fig.~\ref{fig:diferenciaU}),  gauges the truncation error. The control algorithm, by adjusting the size of the timesteps,  keeps this difference around a prefixed value $\tau$, the tolerance.  Hopefully, this tolerance is an indicator of the numerical error. In this paper we define the difference $\mathcal{E}^{(n)}$ in this way:
\begin{equation}
\mathcal{E}^{(n)} = \max_{\text{all}\, k} \left|\widehat{U}_k^{(n)}-{U}_k^{(n)}\right| .
\end{equation}

In the next two subsections we  describe in detail these two methods and discuss their main characteristics and performance. To do this, we  will use the following problem as testbed:
\begin{subequations}
\label{testbed}
\begin{align}
 &\frac{\partial^\gamma u}{\partial t^\gamma}=   \frac{\partial^2 u}{\partial x^2}, \quad 0\le x\le \pi,\\
 &u(x=0,t)=u(x=\pi,t)=0,\\
& u(x,0)=\sin x.
\end{align}
\end{subequations}
Its exact solution is specially simple
\begin{equation}\label{}
u(x,t)= E_\gamma(-t^\gamma) \sin(x)
\end{equation}
with $E_\gamma$ being the Mittag-Leffler function \cite{Mainardi2000}.

\subsection{Trial and error method}

In the trial and error (T\&E) algorithm the procedure for choosing the size of the timesteps is as follows \cite{Quintana-Murillo2013}:
\begin{enumerate}
\item 	 If, initially, $\mathcal{E}^{(n)}$ is larger than the tolerance $\tau$, then we halve the timestep $\Delta_n$  and check whether the new difference $\mathcal{E}^{(n)}$ corresponding to the new timestep (i.e., to the timestep $\Delta_n/2$) is still  larger than the tolerance. We repeat this procedure until the difference $\mathcal{E}^{(n)}$ is smaller than $\tau$. In this case we get the last timestep as the appropriate value.
\item But  if, initially, the difference $\mathcal{E}^{(n)}$ is smaller than the tolerance $\tau$, then we double the size of the timestep. We   repeat this procedure  until $\mathcal{E}^{(n)}$ is larger than the tolerance. When this happens, we take as appropriate timestep the timestep \emph{previous} to the last one.
\end{enumerate}
As starting value for $\Delta_n $ we use the value of the previous timestep $\Delta_{n-1}$. Therefore  $\Delta_0$  is undefined and its value has to be given in order to initiate the algorithm. We have always taken $\Delta_0=0.01$ because this seeding  value is largely irrelevant since the  T\&E algorithm  quickly finds a timestep $\Delta_1$ that leads to a suitable $\mathcal{E}^{(1)}$.

In order to see how good   the T\&E  algorithm is, we  checked its speed (that is, its capacity to integrate the equation over large time intervals employing small CPU times)  and the size of the errors that the method provides when applied to the testbed problem \eqref{testbed}.
Unless otherwise explicitly stated,  the CPU times are not given in seconds but in units of $\mathcal{T}_{50}$, which is the CPU time employed by our method with fixed timesteps (i.e, without implementing any adaptive choice of the size of the timesteps) to get the solution of the problem \eqref{testbed} for $\gamma<1$ when $50$ timesteps are used (in our computer $\mathcal{T}_{50} \approx 1.4 $ seconds). In this way the CPU-time values reported here are expected to be roughly independent of the particular computer system employed.   The normalized CPU time required to evaluate the numerical solution of a given problem  at time $t$  will be denoted by $T_\text{CPU}(t)$.

\begin{figure}[th]
\begin{center}
        \includegraphics[width=0.42\textwidth,angle=0]{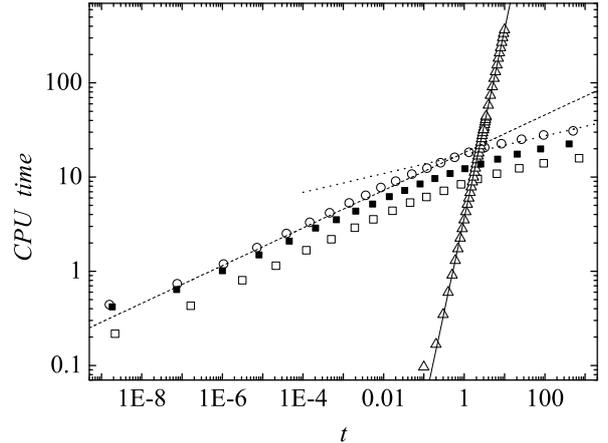}
\end{center}
\caption{\label{CPUvsTimeFijovsTEvsPred}   Normalized computational time $T_\text{CPU}(t)$ required by  the fixed-step method  with $\Delta=0.01$ (triangles), by the T\&E method (circles), and by the predictive method with $\theta=3/2$  and $\omega=1$ (open squares) and $\omega=1/2$ (solid squares),  to solve problem \eqref{testbed} with $\gamma=1/4$ up to time $t$. In all cases $\Delta x=\pi/40$  and $\tau=10^{-4}$.  The lines  are guides to the eye; their slopes (0.2 for the dashed line, 0.1 for the dotted line, 2 for the solid line)  provide estimates of the power exponent $\beta$ in $T_\text{CPU}(t) \sim t^\beta$.  }
\end{figure}

In Fig.~\ref{CPUvsTimeFijovsTEvsPred}, we compare   $T_\text{CPU}(t)$  for $\gamma=1/4$ corresponding to the T\&E method with tolerance $\tau=10^{-4}$   with the  $T_\text{CPU} $ values for the standard (non-adaptive) method with fixed timesteps of size $\Delta_n=0.01$.  First we see that, as expected, the CPU time required by the standard method grows quadratically: $T_\text{CPU}\propto t^2$.   However, for the T\&E  method one finds that   $T_\text{CPU}\propto t^\beta$ with $\beta\approx 0.2$. That is,  the growth of the computational time is not quadratic, not even linear, but  sublinear!  In fact,  the growth is even slower for longer times ($\beta \approx 0.1$)! This means that, except for short times, the adaptive method is immensely faster than the standard method with fixed timesteps.  In Table \ref{tablaCPU} we give some specific values of $T_\text{CPU}(t)$ in seconds. It is clear that the standard method becomes useless for times $t$ above a few tens when $\Delta_n=0.01$ or, equivalently, when the number of timesteps is above a few thousands.

\begin{table}
\centering
\begin{tabular}{|c|c|c|c|c|c|c|}
  \hline
  $t$ & 0.12 & 1.29& 2.58 &3.22 &  8.34 & 502.4\\
 $T_\text{CPU}$ (FT)   &   0.2      &8  & 33& 49   & 346 & $\gtrsim$ 2 weeks \\
   $T_\text{CPU}$  (T\&E)  & 17   &26  &29 &29 &   32 & 43  \\
  \hline
\end{tabular}
\caption{CPU time (in \emph{seconds}) employed to get the solution of \eqref{testbed} with  $\gamma=1/4$ up to time $t$ by (i) the method with fixed timesteps (FT) with $\Delta=0.01$, and (ii) by the T\&E method  with tolerance $\tau=10^{-4}$. In both methods $\Delta x=\pi/40$. The CPU time required by the FT method to find the solution at $t=502.4$ is an estimate obtained by extrapolation (see Fig.~\ref{CPUvsTimeFijovsTEvsPred}). }
\label{tablaCPU}
\end{table}

\begin{figure}[th]
\begin{center}
\includegraphics[width=0.42\textwidth,angle=0]{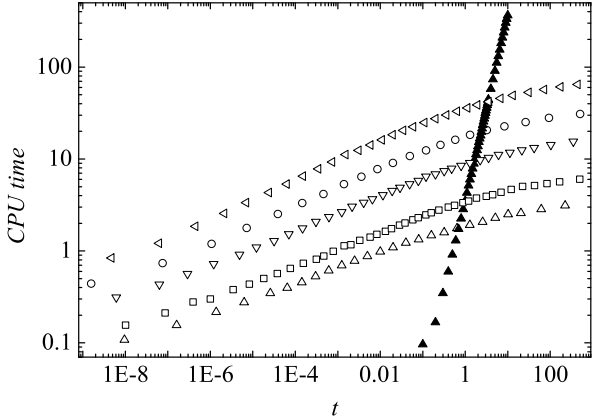}
\end{center}
\caption{\label{CPUvsTimeFijovsTEvarias}   Normalized computational time $T_\text{CPU}$ vs $t$ for the method with fixed timesteps  with $\Delta=0.01$ (solid triangles) and for the T\&E method  with tolerance $10^{-5}$ (left triangles), $5\times 10^{-4}$ (circles), $10^{-4}$ (down triangles), $2\times 10^{-4}$ (squares), $10^{-3}$ (open up triangles).
In all cases $\gamma=1/4$ and  $\Delta x=\pi/40$.
}
\end{figure}

The CPU times of the T\&E method given in Fig.~\ref{CPUvsTimeFijovsTEvsPred} and Table \ref{tablaCPU} correspond to a tolerance $\tau=10^{4}$.   In Fig.~\ref{CPUvsTimeFijovsTEvarias},  we show   $T_\text{CPU}(t)$  for several values of $\tau$. As expected,  the speed of the method increases when the tolerance decreases (of course, the price to be paid is that  the method is then less accurate as we shall see below).  It is interesting to note that Fig.~\ref{CPUvsTimeFijovsTEvarias}  shows that the CPU time is, in fact, roughly proportional to $\tau^{-\eta}$ with $\eta$ around 1. For example, when $\tau$ changes from $\tau=10^{-3}$ to $\tau=10^{-4}$, ones sees that the CPU time increases approximately by a factor of ten.   For other values of $\gamma$, a very similar behaviour is found. This can be explained by the following  back-of-the-envelope   argument. By construction, $\tau\approx  \mathcal{E}^{(n)}$, but, as  will be seen in Section \ref{predictiveSec}, $\mathcal{E}^{(n)} \sim \Delta_n ^\theta$ where $\theta\approx 3/2$. Then $ \Delta_n   \sim \tau^{1/\theta}$. For a given time, $t$, $T_{CPU}(t)\sim n^2\sim (t/\Delta t)^2$ where $\Delta t$ is here the average value of the timesteps given until time $t$. But if the size of the timesteps scales roughly as $\tau^{1/\theta}$, then one expects that its average should scale similarly, $\Delta t \sim  \tau^{1/\theta}$. Therefore $T_{CPU}(t)\sim \tau^{-\eta}$ with $\eta=2/\theta$. For $ \theta \approx  3/2$ one gets $\eta\approx 1.3$, which is not far from 1.

\begin{figure}[h]
\begin{center}
        \includegraphics[width=0.42\textwidth,angle=0]{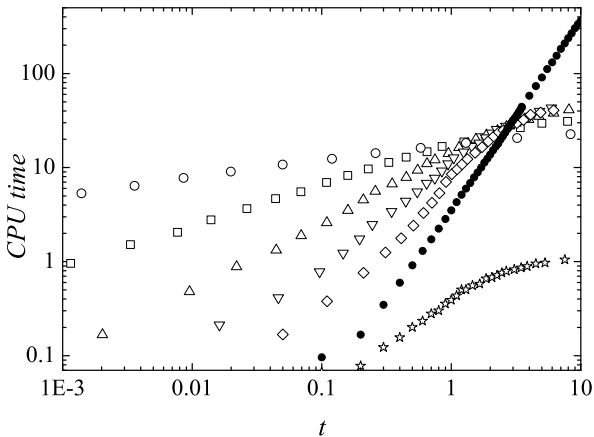}
\end{center}
\caption{\label{CPUvsTimevsGammaTE}  Normalized computational time $T_\text{CPU}$ vs $t$ when problem \eqref{testbed} is solved by means of the T\&E method with $\tau=10^{-4}$ for $\gamma=0.25, 0.5, 0.75, 0.9, 0.99, 1$  (open circles, squares, up triangles, down triangles, diamonds, stars, respectively) and by the method with fixed timesteps  with $\Delta=0.01$ (solid circles). In all cases $\Delta x=\pi/40$.
}
\end{figure}

Figure~\ref{CPUvsTimevsGammaTE} shows how the computational time $T_\text{CPU}(t)$ depends on the fractional parameter $\gamma$.  One sees that,  approximately, $T_\text{CPU}(t) \sim t^\beta$ and that $\beta$ increases when $\gamma$ increases ($\beta\to 2$ when $\gamma\to 1^-$). However,  the value $\gamma=1$ is singular: there is a  drastic change of the CPU times between $\gamma \to 1^-$ and $\gamma= 1$. The reason for this is clear:  for $\gamma=1$ the fractional derivative  is a pure differential operator (i.e., a local operator, no longer an  integro-differential operator), and then it is not necessary to carry out the sum of the right hand side of Eq. \eqref{ImpliAlg1}, which is what makes fractional finite difference methods so (increasingly) slow.   For this reason, the CPU times required to integrate the normal diffusion problem ($\gamma=1$) are far shorter than the CPU times for subdiffusion problems ($\gamma<1$).   Note that for $\gamma<1$, no matter how close is $\gamma$ to unity, one has to spend a lot of computation time evaluating this sum even though one knows that the closer  $\gamma$ is to unity the closer  the sum is to zero.

\begin{figure}[h]
\begin{center}
        \includegraphics[width=0.42\textwidth,angle=0]{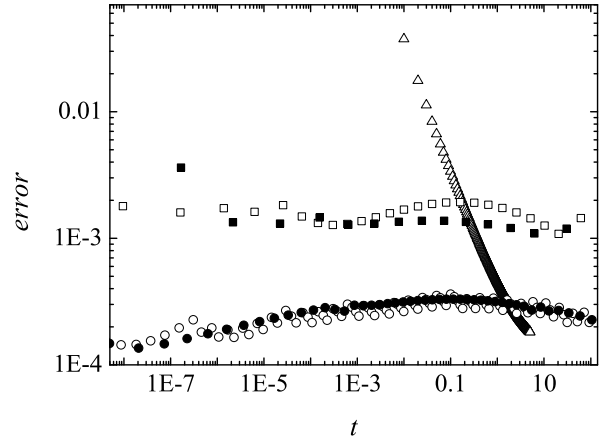}
\end{center}
\caption{\label{ErrorvsTime1}  Numerical error   vs.~$t$ when problem \eqref{testbed} is solved by means of the fixed step method with $\Delta_n  =0.01$  (triangles),   the T\&E method with tolerance $10^{-4}$ (circles) and $10^{-3}$ (squares), and   the predictive method with  $\omega=1/2$, and tolerance $10^{-4}$ (solid circles) and $10^{-3}$ (solid squares). In all cases $\Delta x=\pi/40$.  }
\end{figure}

Regarding the  error,  one sees in Fig.~\ref{ErrorvsTime1}  that the adaptive algorithm provides quite homogeneous errors, that is, this technique has the convenient property that it keeps the errors to a desired degree of accuracy, neither too large nor too small.  This should be compared with the quite uneven errors of the standard method with fixed timesteps.
(In short, non-uniform timesteps lead to uniform errors while uniform  timesteps lead to quite non-uniform errors.)
Finally,  one sees in Fig.~\ref{ErrorvsTime1}  that the errors are close to the tolerance, in particular, that they are of the same order of magnitude (around three times the tolerance in this case).  This example illustrates the  fact that this adaptive algorithm  has the nice property that the tolerance, a quantity one can fix at will,  is a convenient \emph{indicator} of the accuracy of the numerical method.

Indeed, this method is both fast and accurate. In the next subsection, we present another adaptive method, the predictive method,  that is as accurate as the T\&E  method but a bit faster.

\subsection{Predictive method}
\label{predictiveSec}

The adaptive predictive method is also based on the step-doubling technique.   The starting point is to assume that the difference   $\mathcal{E}^{(n)}$ scales as a power of the size of the timesteps:
\begin{equation}
\label{Endt}
\mathcal{E}^{(n)}\sim  \Delta_n^\theta .
\end{equation}
Provided  this relationship holds, and from the value of the difference $\mathcal{E}_\text{old}^{(n)}$ for a given timestep $\Delta_n^\text{old}  $, one can easily predict the size of the timestep $\Delta_n^\text{pred}  $ that leads to an error equal to the tolerance, $\mathcal{E}_\text{pred}^{(n)}= \tau$ , namely,
\begin{equation}\label{dtnew}
 \Delta_n^\text{pred}  = {\Delta_n^\text{old}}    \,
   \left[ \frac{ \tau}{ \mathcal{E}_\text{old}^{(n)} }\right]^{1/\theta}  .
\end{equation}
One expects that the direct use of this timestep should spare one from wasting computer time trying to find the right timestep (the one that leads to a difference $\mathcal{E}^{(n)}$  of the order of the tolerance) by means of a blind succession of trials and errors of the size of the timestep as the T\&E method does.  This prompts us to propose the following predictive step-doubling algorithm:
\begin{enumerate}
\item   If, for the initial  value of $\Delta_n$, the difference satisfies
\begin{equation}\label{predcon}
   \tau/2\le \mathcal{E}_n \le 2 \tau,
\end{equation}
then this timestep is accepted.

\item Otherwise, we employ a new timestep given by the formula
\begin{equation}\label{predDelta}
 \Delta_n^\text{new}  =  \omega\,  {\Delta_n^\text{old}}   \,
   \left[ \frac{ \tau}{ \mathcal{E}_\text{old}^{(n)} }\right]^{1/\theta} + (1-\omega)\,\Delta_n^\text{old}
\end{equation}
until the corresponding difference $\mathcal{E}_\text{new}^{(n)}$ satisfies condition \eqref{predcon}.
\end{enumerate}
The starting value for $\Delta_n $ is just  $\Delta_{n-1}$, and we take $\Delta_0=0.01$ in all cases.   As mentioned previously, the particular initial value for $\Delta_0$  is largely irrelevant because the above step-doubling algorithm,  as does the T\&E algorithm, quickly finds the right timestep $\Delta_0$. The parameter $0<\omega\le 1$ is a kind of under-relaxation parameter \cite{Press:2007:NRE:1403886,LeVeque2007}  that affects the speed and robustness of the algorithm. We have found that the (pure predictive) parameter $\omega=1$ usually leads  to the fastest adaptive algorithm (see Fig.~\ref{CPUvsTimeFijovsTEvsPred}). However, in some extreme cases,  the pure predictive method breaks down  because the choice of the timesteps enters into an infinite loop   $\Delta_n^a \to  \Delta_n^b    \to \Delta_n^a   \ldots$ due to the fact that their  corresponding  differences $\mathcal{E}^{(n)}$  never fulfill  the (exit) condition \eqref{predcon}.
When this happens, we have found that the use of an under-relaxation value ($0<\omega<1$) fixes this problem. In particular, we have never found this kind of problem for  $\omega=1/2$.

\begin{figure}[th]
\begin{center}
        \includegraphics[width=0.42\textwidth,angle=0]{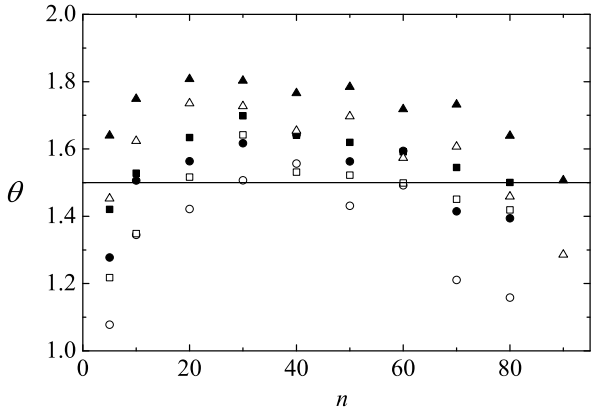}
\end{center}
\caption{\label{Thetavsn}  Scaling exponent $\theta$ of the predictive method vs.~the number of timesteps when problem \eqref{testbed} is solved for several values of $\gamma$ and  $\Delta_n=m \Delta_{n-1}$   (open symbols) and  $\Delta_n=m \Delta_{n-1}/3$ (solid symbols) with $m=1,2,\ldots 10$. The line marks the value  $\theta=3/2$.
   In all cases $\Delta x=\pi/80$ and $\omega=1$.}
\end{figure}

\begin{figure}[th]
\begin{center}
        \includegraphics[width=0.42\textwidth,angle=0]{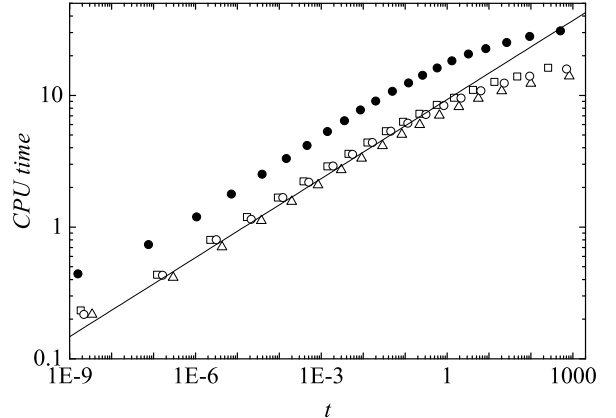}
\end{center}
\caption{\label{CPUvsTimeTEvsPred}  Normalized computational time $T_\text{CPU}$ vs.~$t$ when problem \eqref{testbed} is solved by means of the T\&E method (stars) and  the predictive method with   $\theta=3/2$ (squares), $\theta=5/4$ (up triangles), and  $\theta=5/3$ (circles).  In all cases $\Delta x=\pi/40$ and $\omega=1$. The line, of slope 0.2,  is a guide to the eye.
  }
\end{figure}

The predictive method relies  on Eq.~\eqref{Endt} so that it is pertinent to check the validity of this power-law scaling.
Figure \ref{Thetavsn}    shows the values of $\theta$ obtained by fitting $\mathcal{E}^{(n)}$ to Eq. \eqref{Endt} for several values of $\gamma$,   number of timesteps $n$, and values of $\Delta_n$. Specifically, we  studied how $\mathcal{E}^{(n)}$ scales with the size of the the timestep employing for $\Delta_{n}$ the values $m \Delta_{n-1}$  and also  $m \Delta_{n-1}/3$  with $m=1,2,\ldots 10$, with  $\Delta_{n-1} $ being the size of the previous timestep.  We found the exponent $\theta$ is always between one and two.  As a simple overall effective value, we  will always use $\theta=3/2$ in this paper.  Of course,  this choice would be questionable if the method were very sensitive to the specific value of $\theta$.  It turns out that this is not the case. Regarding the computational time, Fig.~\ref{CPUvsTimeTEvsPred} shows that the specific value of $\theta$ is hardly relevant. Note, however,   that the predictive method is  \emph{faster} than the T\&E method, around three times faster in this example.
Regarding the errors, one sees in  Fig.~\ref{ErrorvsTime1} that their behaviour is  quite similar to the behaviour of the errors of the T\&E method. For both methods the errors are nicely close to the prefixed tolerance.

\section{Two examples. Source term and non-homogeneous boundary conditions}
\label{twoExamples}

\subsection{A problem with a steep source term }
A relevant feature of the  two adaptive methods we have presented above is that they can dynamically adapt the size of the timesteps according to the behaviour of the solution. For example, if at a given time we introduce an external perturbation into the system, we expect the adaptive method to be able to take care of this through the use of a temporal mesh suited to the behaviour of the perturbation.  In this subsection we give a clear example of this. The problem we consider is  Eq.~\eqref{ecu1} with $K=1$,  boundary conditions $u(0,t)=u(\pi,t)=0$,  and the source term
\begin{equation}\label{fxt}
 f(x,t)=\left[1+\frac{\Gamma(1+p)t^{-\gamma}}{\Gamma(1+p-\gamma)}  \right] a  t^{p} \sin x .
\end{equation}
Its exact solution is
\begin{equation}
u(x,t)=\left[E_\gamma(-t^\gamma)+a t^p\right]\, \sin x  .
\end{equation}
\begin{figure}[th]
\begin{center}
        \includegraphics[width=0.42\textwidth,angle=0]{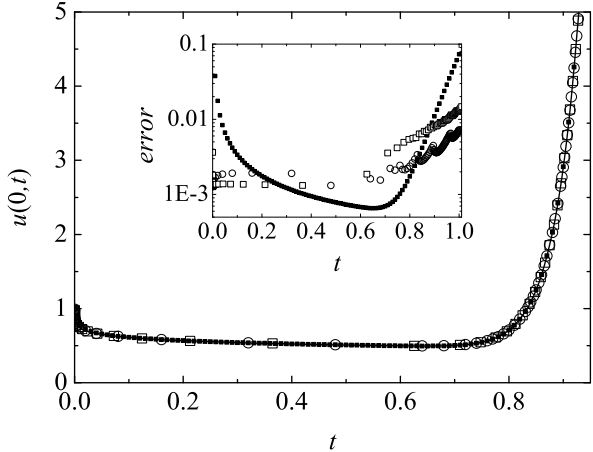}
\end{center}
\caption{\label{casoFuente}  Solution and numerical errors  at the midpoint $u(\pi/2,t)$ of the problem described in the main text for $\gamma=1/4$ and source term \eqref{fxt} with $a=p=20$.  Solid squares: numerical method with $\Delta_n=0.01$; circles: T\&E method with tolerance $\tau=10^{-3}$; open squares: predictive method with  $\tau=10^{-3}$ and $\omega=1/2$; line: exact solution.  In all cases $\Delta x=\pi/40$.
 }
\end{figure}
Here we take $a=p=20$ because this leads to a solution with three different time regimes: the short-time regime where the solution changes very fast, the intermediate regime (roughly until times a little below $t=1$) where the solution changes very slowly, and a final regime for longer times where again the solution changes very fast.  This case is then a good example with which to test adaptive methods. The results provided by the T\&E method, the predictive method with $\omega=1/2$, and the method with fixed timesteps are shown in Fig.~\ref{casoFuente}. Remarkably, we have found that, in some cases, the pure predictive method ($\omega=1$) breaks down in this example as described in Section \ref{predictiveSec}.  It is also remarkable the way in which the size of the timesteps of the  adaptive methods changes  according to the behaviour of the solution: for small times and for times around $t=1$, the solution change very fast and then the adaptive methods tiptoe  in these regions keeping the numerical errors small; for intermediate times the solution changes very slowly, and the adaptive methods react by making large strides, thus going fast forward in time  although not at the expense of increasing the numerical errors.

\subsection{A problem with non-homogeneous boundary conditions }

We want to find the density profile at any time of a set of continuous-time random walkers moving in a  one-dimensional finite medium, initially void of walkers,  when there is a reservoir of walkers at one end of the medium (so that their concentration is constant there)  and they are  completely removed from the system at the other end.
In mathematical terms, the problem we have to solve then is given by Eq.~\eqref{ecu1} with $0\le x\le L$, boundary conditions $u(0,t)= u_0$, $u(L,t)=0$, and initial condition $u(x,0)=0$.
The exact solution can be obtained by solving the problem in the Laplace space,  or directly by means of the method of images \cite{Metzler2000}:
\begin{align}
u(x,t)&= u_0  \sum_{m=0}^M  H_{10}^{11} \left[m z_c+ z\LM{1,\gamma/2}{0,1}\right]\nonumber\\
&-u_0 \sum_{m=1}^M  H_{10}^{11} \left[m z_c- z\LM{1,\gamma/2}{0,1}\right], \quad M\to\infty,
\label{soluEscalon}
\end{align}
where $z = x/(K t^{\gamma})^{1/2}$, $z_c = 2 L/(K  t^{\gamma})^{1/2}$, and $H_{10}^{11}$ is a Fox $H$ function  \cite{Metzler2000,Mathai1978}.  When $\gamma=1$, the Fox function becomes the complementary error function $H_{10}^{11}\left[z\left|^{1,1/2}_{0,1}\right.\right]=\text{erfc}(z/2)$, and the classical solution \cite[Eq.~(6), p.~310]{Carslaw1959} is recovered.

\begin{figure}[th]
\begin{center}
        \includegraphics[width=0.42\textwidth,angle=0]{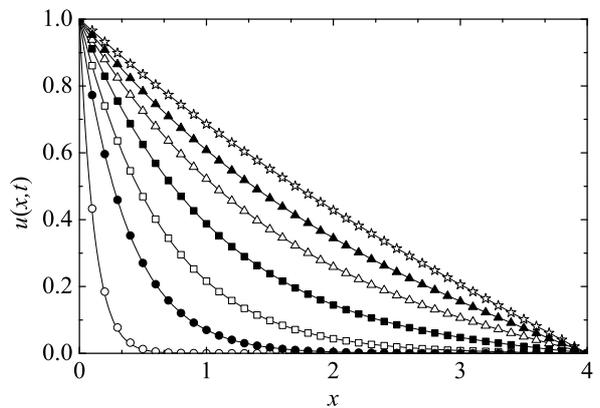}
\end{center}
\caption{\label{casoEscalon}  Solution $u(x,t)$ vs.~$x$ of the problem with nonhomogeneous boundary condition described in the main text for  $\gamma=1/4$, $K=1$, $L=4$, and several values of $t$. The lines are the exact solution given by Eq.~\eqref{soluEscalon} with $M=8$, the open symbols are the numerical solutions obtained by means of the T\&E method, and the solid symbols correspond to the predictive method with  $\omega=1/2$. For both methods $\tau=10^{-3}$.  The solutions correspond to  (from left to right)  $t=1.91\times 10^{-8}, 2.67 \times 10^{ -4}, 2.00 \times 10^{ -2}, 8.93 \times 10^{ -1}, 2.05\times 10^{1}, 2.68 \times 10^{2 }, 1.14 \times 10^{4 }$.  In all cases,  $\Delta x=0.1$ and $\tau=10^{-3}$.
 }
\end{figure}

In Fig.~\ref{casoEscalon}, we compare the exact solution with the numerical solution obtained with the T\&E and predictive methods.
This problem makes clear the immense advantage of adaptive methods over methods with fixed timesteps.  In order to find the solution close to the stationary state, one has to reach times around $10^4$. Therefore, in order to get this solution  by means of a reasonable number  of timesteps of fixed size, one has to use large timesteps,  let us say $10^4$ timesteps of size $\Delta_n =1$, which means that all the changes of the system from $t=0$ up to $t=\Delta_n =1$ would be overlooked.  In our case, see Fig.~\ref{casoEscalon}, this would mean overlooking a time interval in which substantial and relevant changes in the solution occurs.
In other words,  Fig.~\ref{casoEscalon} shows that, in order to conveniently track the solution from the initial condition to the stationary state, one has to employ times that span twelve orders of magnitude (from $t\sim 10^{-8}$ to $t\sim 10^4$). No computer employing finite difference methods with fixed timesteps can handle this problem in a reasonable computation time.

\section{Summary and conclusions}

A mayor drawback of finite difference methods for fractional equations is that obtaining the solution for every new timestep is increasingly costly in terms of computational time. This implies that the number of timesteps required to find a solution should be as small as posible while keeping a reasonable accuracy in the numerical solution.  In this respect, finite difference methods with homogeneous timesteps have the additional drawback that  they typically provide numerical results with  quite uneven accuracy (small in some time regions but relatively large in others).

In order to lessen these two problems, we have proposed the use of adaptive methods with adaptive timesteps. This kind of method has the great advantage that the size of the timesteps can be tailored  to the behaviour of the solution. For example, one can choose small timesteps  only   when the solution is changing rapidly in order to keep the accuracy of the method  and  to  track down the relevant features of the solution at these time scales. On the other hand,  one can choose large timesteps when the solution changes slowly.  In this way, without losing accuracy, the method can advance  with large strides when feasible so that long times can be reached.

We have studied two different adaptive methods. Both employ an integration algorithm based on the L1 discretization of the Caputo fractional derivative with non-homogeneous timesteps, but differ in their adaptive algorithm, i.e., in the way in which the size of every timestep is determined. Both methods are fast (immensely faster than the corresponding non-adaptive method with fixed timesteps in many cases) and provide solutions with an accuracy that, to a large extent, is consistent throughout the integration time interval. Remarkably, this accuracy can be pre-set at will through setting up a kind of self-consistent parameter (the tolerance) that is in many cases an excellent indicator of the final accuracy of the numerical solution.

\section*{Acknowledgments}

This work was partially funded by the Ministerio de Ciencia y Tecnolog\'ia (Spain) through Grant No. FIS2010-16587 (partially financed by FEDER funds) and by the Junta de Extremadura through Grant. No. GRU10158.


\end{document}